\documentclass[10pt, twoside]{amsart}
\usepackage{graphicx}
\usepackage{caption}
\usepackage{subcaption}
\usepackage[utf8]{inputenc}
\usepackage{amssymb} 
\usepackage{amsmath}
\usepackage{algorithm}
\usepackage{algpseudocode}
\usepackage[noadjust]{cite}    
\usepackage{epsfig}

\newtheorem{theorem}{Theorem}[section]

\theoremstyle{definition}

\newtheorem{remark}[theorem]{Remark}

\numberwithin{equation}{section}

\newcommand{\dt}{\;{\rm d}t}

\begin{document}

\title[Deep learning using the Brezis--Ekeland principle]{
Deep learning for gradient flows using the Brezis--Ekeland principle
}

\author[L. Carini]{Laura Carini}
\address{Dipartimento di Mathematica, Universit\`a di Trento,
38123 Trento, Italy} 
\email{laura.carini@studenti.unitn.it}

\author[M. Jensen]{Max Jensen}
\address{Mathematics Department, University College London, 25 Gordon Street, London, WC1H 0AY, United Kingdom}
\email{max.jensen@ucl.ac.uk}

\author[R. N\"urnberg]{Robert N\"urnberg}
\address{Dipartimento di Mathematica, Universit\`a di Trento,
38123 Trento, Italy} 
\email{robert.nurnberg@unitn.it}

\thanks{}

\begin{abstract}
We propose a deep learning method for the numerical solution of partial differential equations that arise as gradient flows. The method relies on the Brezis--Ekeland principle, which naturally defines an objective function to be minimized, and so is ideally suited for a machine learning approach using deep neural networks. We describe our approach in a general framework and illustrate the method with the help of an example implementation for the heat equation in space dimensions two to seven.
\end{abstract}



\maketitle

\section{Introduction} \label{sec:1}

In this paper we advocate a deep learning approach for solving parabolic partial differential equations (PDEs) 
\[
u_t + \partial \phi(u) = f,
\] 
that arise as evolution equations for gradient flows. We exploit the variational principle of the seminal papers by Brezis and Ekeland, \cite{BrezisE76,BrezisE76a}, now commonly known as the Brezis--Ekeland principle, \cite{Stefanelli08,Roubicek13}. 

Using neural networks for the numerical solution of PDEs has become increasingly popular over the last decade. An attraction of neural network-based approaches is their suitability for high-dimensional problems. 
For a comprehensive review of current developments we refer to \cite{Karniadakis2021422, BlechschmidtE21}. Among currently popular techniques are methods based on residual minimization, e.g.~see  \cite{SirignanoS18,RaissiK18,RaissiPK19}, and on the reformulation as a backward stochastic differential equation, e.g.~\cite{Han:2016wf, Han:2018kv, HenryLabordere:2017fe, Raissi:2018wq}.

Most relevant for this work are variational approaches. For an elliptic problem, E and Yu, \cite{EY18}, proposed a deep learning method based on a variational principle that leads to a natural optimization framework. 
An approach connecting variational principles with convex duality for stationary equations was taken in \cite{KaltenbachZ22preprint}, which mirrors some aspects of our work for gradient flows.

The aim of our approach is three-fold:
\begin{enumerate}
    \item Machine learning approaches have been criticized for their less developed methodology to bound, or at least estimate, the approximation error. It is therefore interesting that the minimum of the Brezis--Ekeland functional, which is minimized during the learning process, is guaranteed to be zero for the exact solution, thus providing an error measure that is known at the point of computing the neural net approximation.
    \item Adversarial networks have very successfully been applied across multiple problem classes of machine learning. We were intrigued by the question whether duality can be a context in which the concept of adversarial networks is translated to partial differential equations as well as to convex analysis, by introducing a neural network for the primal and another one for the dual problem. In a resulting min-max formulation the training stages of the respective networks take opposing, or adversarial, roles in finding the value of the joint loss functional.
    \item Finally, we wish to construct a method which takes advantage of the specific structural properties of gradient flows, based on the relevance of these properties in the literature for the construction of finite element methods for time-dependent PDEs.
\end{enumerate}

The outline of the remainder of the paper is as follows: in section 2 we introduce gradient flows and the Brezis--Ekeland principle; in section 3 we formulate our deep learning approach; in section 4 we discuss the computer implementation of the method; in section 5 we present numerical experiments, followed by conclusions.

\section{The Brezis--Ekeland principle for gradient flows} \label{sec:BEP}
Let $V \subset H \subset V^*$
be a Gelfand triple and $T>0$ a fixed time. In addition, let $\phi\in C^1(V)$ be convex, $f \in L^2(0,T;V^*)$ and $u_0 \in V$. We consider the gradient flow
\begin{equation} \label{eq:gradflow}
    u_t + \partial \phi(u) = f \quad \text{a.e.\ in }\ (0,T), \quad
    u(0) = u_0.
\end{equation}
The Brezis--Ekeland principle asserts that solutions $u \in Y$ to \eqref{eq:gradflow}, with
\[
    Y := \{ w \in L^2(0,T;V) \cap H^1(0,T;V^*) \; : \, w(0) = u_0 \},
\]
are the global minimizers of the functional $\Phi : Y \to [0,\infty]$ defined by
\begin{align} \label{eq:Phi}
    \Phi(w) = & \tfrac12 \|w(T)\|^2_H - \tfrac12 \|w(0)\|^2_H + \int_0^T \phi(w) + \phi^*(f-w_t) - \langle f,w \rangle \dt \nonumber \\
            = & \int_0^T \phi(w) + \phi^*(f-w_t) + \langle w_t - f,w \rangle \dt,
\end{align}
using $\int_0^T \langle w_t,w \rangle \dt = \tfrac12 \|w(T)\|^2_H - \tfrac12 \|w(0)\|^2_H$. Here $\langle\cdot,\cdot\rangle$ is the duality pairing between $V$ and $V^*$, and $\phi^*(w) = \sup_{v \in V} \langle w,v \rangle - \phi(v)$ is the conjugate of $\phi$. In fact, owing to \cite[Theorem 8.99]{Roubicek13}, $u$ solves \eqref{eq:gradflow} if and only if 
\begin{equation} \label{eq:BEPorig}
    \Phi(u) = \min \{ \Phi(w) \; : \, w \in Y \} = 0.
\end{equation}
Inspired by the work in \cite{Stefanelli09}, we now consider a time-discrete variant of \eqref{eq:Phi} and the associated minimization problem. For that purpose, we divide $[0,T]$ into $N$ sub-intervals with end points $t_0 = 0 < t_1 < \cdots < t_N = T$. 

The parabolic nature of \eqref{eq:gradflow} ensures that $u(t)$ only depends on $u(s)$ if $s \le t$, but not if $s > t$. Together with the Brezis--Ekeland principle \eqref{eq:BEPorig} this guarantees that minimization and summation may be interchanged:
\begin{align} \label{eq:BEP} 
         & \min_{w\in Y} \Phi(w)\\
    = \, & \min_{w\in Y} \sum_{n=1}^N \int_{t_{n-1}}^{t_n} \phi(w) + \phi^*(f-w_t) + \langle w_t - f,w \rangle \dt \nonumber \\ 
    = \, & \sum_{n=1}^N \min_{\substack{w_n\in Y\\ w_n(t_{n-1}) = w_{n-1}(t_{n-1})}} \int_{t_{n-1}}^{t_n} \phi(w_n) + \phi^*(f-(w_n)_t) + \langle (w_n)_t - f,w_n\rangle \dt, \nonumber
\end{align}
where we have defined $w_0(0) = u_0$. The representation \eqref{eq:BEP} of the minimization problem suggests a clear strategy for our deep learning method: we will sequentially solve $N$ optimization problems for the PDE \eqref{eq:gradflow} on the time intervals $[t_{n-1}, t_n]$, where the initial data is either given by $u_0$ at the first step, or by the previously computed solution at time $t_{n-1}$.

\subsection*{The heat equation}

The canonical example of a gradient flow is the heat equation. Given a domain $\Omega \subset \mathbb{R}^d$, $d \geq 1$, we consider the PDE:
\begin{equation}\label{eq:PDEpbl}
    \begin{cases}
        u_t - \kappa\Delta u = f, \quad & \text{in } (0,T) \times \Omega,\\
        u(0,\cdot) = u_0, \quad & \text{in } \Omega,\\
        u = 0 , \quad & \text{on } (0, T) \times \partial\Omega,
    \end{cases}
\end{equation}
where $\kappa > 0$. Upon defining $H=L^2(\Omega)$, $V = H^1_0(\Omega)$ and $\phi(u) = \frac\kappa2 \|\nabla u\|^2_{L^2}$, the problem (\ref{eq:PDEpbl}) is a special case of \eqref{eq:gradflow}, and the Brezis--Ekeland functional  \eqref{eq:Phi} in this case reduces to:
\begin{equation}\label{eq:BEphi1}
    \Phi(w) = \int_0 ^T \frac\kappa2 \|\nabla w\|^2_{L^2} +  \frac1{2\kappa} \left[ \sup_{v \neq 0}\frac{ \langle f-w_t, v \rangle }{  \| \nabla v  \|_{ L^2}} \right]^2 + \langle w_t -f , w \rangle \dt,
\end{equation}
where we have used that $\phi^*(w) = \frac1{2\kappa} \left[ \sup_{v \in H^1_0(\Omega) \setminus \{0\}}\frac{ \langle w, v \rangle }{  \| \nabla v  \|_{ L^2}} \right]^2$, e.g.~see \cite{Roubicek13}. We remark that $\phi^*$ defines a norm on $H^{-1}(\Omega) = (H^{1}_0(\Omega))^*$. Combining \eqref{eq:BEP} and \eqref{eq:BEphi1} we find that $\min_{w\in Y} \Phi(w)$ equals
\begin{align} \label{eq:BEPheat}
    \sum_{n=1}^N \min_{\substack{w_n(t_{n-1})\\ = w_{n-1}(t_{n-1})}}
    \int_{t_{n-1}}^{t_n} \frac\kappa2 \|\nabla w_n \|^2_{L^2} + \frac1{2\kappa} \left[ \sup_{v \neq 0}\frac{ \langle f-(w_n)_t, v \rangle }{  \| \nabla v  \|_{ L^2}} \right]^2 + \langle (w_n)_t - f,w_n\rangle \dt, 
\end{align}
with the solution $u$ of \eqref{eq:PDEpbl} being a minimizer in the sense that $\Phi(u) = 0$ and that the choice $w_n = u$, $n \in \{1, \ldots, N \}$, yields a minimizer of \eqref{eq:BEPheat} over~$Y^N$.

\section{The deep learning approach}

In this section we discuss a deep learning algorithm to find approximations of the solution to \eqref{eq:gradflow}. For simplicity we restrict our attention to the heat equation, so that \eqref{eq:BEPheat} is our starting point.

We wish to find the approximations $u^n_h \approx u(t_n, \cdot)$, $n=1,\ldots,N$, where $u^n_h$ is given by a neural network. More generally, in this paper we view a neural network $\hat{u}_h$ as a function determined through its weights $\theta$. Given $\theta$ the neural network takes a position $x \in \Omega$ as input and returns $\hat{u}_h(x; \theta)$ as output. The approximation set containing the functions $u^n_h$ is thus given by
\[
    \mathbb U_h := \{ \hat{u}_h(\cdot; \theta) \; : \, \theta \in \Theta \},
\]
where $\Theta$ is the set of possible weights. In this notation $u^n_h = \hat{u}_h(\cdot; \theta)$ is the neural network with a choice of weights $\theta$ determined through the method described in this section.

In order to define a discrete version of the Brezis--Ekeland functional \eqref{eq:BEPheat} and perform its minimization over $\mathbb U_h$, we turn our attention to the interpretation of $v$ in $\phi^*(w) = \frac1{2\kappa} \left[ \sup_{v \in H^1_0(\Omega) \setminus \{0\}} \langle w, v \rangle/\| \nabla v  \|_{ L^2} \right]^2$. Also the $v$ are approximated by neural networks. Since their architecture may be different compared to $\hat{u}_h$, we introduce $\hat{v}_h(x; \eta)$ and
\[
    \mathbb V_h := \{ \hat{v}_h(\cdot; \eta) \; : \, \eta \in \mathrm{H} \},
\]
where $\mathrm{H}$ is the set of possible weights in $\hat{v}_h$. 

While we assume $\mathbb U_h, \mathbb V_h \subset H^1(\Omega)$ throughout, elements of $\mathbb U_h$ and $\mathbb V_h$ will in general not belong to $H^1_0(\Omega)$ because the Dirichlet boundary conditions may not be satisfied homogeneously. Therefore, similarly to \cite{EY18}, we introduce a penalty term into $\phi^*(w)$ to obtain the functional
\begin{equation} \label{eq:phi_h}
    \phi^*_h(w_h) = \frac{1}{2 \kappa} \left[\sup_{v_h \in \mathbb V_h \setminus \{0\}}
\frac{ (  w_h , v_h ) }{ \left( \| \nabla v_h  \|_{ L^2}^2 + \lambda \|v_h\|^2_{L^2(\partial\Omega)} \right)^{\frac12}} \right]^2,
\end{equation}
where $\lambda > 0$ is a penalty parameter depending on the structure of the neural net. Then the denominator in \eqref{eq:phi_h} cannot vanish for $v_h \neq 0$ due to a Poincar\'e-Friedrichs inequality and the penalization weakly imposes homogeneous boundary conditions on any maximising $v_h$ as $\lambda \to \infty$.

It remains to discretize the time derivatives in the Brezis--Ekeland functional. To this end we substitute $(w_n)_t$ in \eqref{eq:BEPheat} by backward time differences. Let $\Delta t_n = t_n - t_{n-1}$ and let $(\cdot,\cdot)$ denote the $L^2$-inner product over $\Omega$. Inspired by \eqref{eq:BEPheat}, we then define the solution of the deep learning method through the following sequence of optimization problems: Given $u_h^{n-1} \in \mathbb U_h$, for $n =1,\ldots,N$, find a minimizer
$u_h^n \in \mathbb U_h$ to
\begin{align} \label{eq:Phi_n}
    \Phi_n(w_h) = \, &  \frac{\kappa\Delta t_n}2 \|\nabla w_h\|^2_{L^2} + \Delta t_n \phi^*_h\Bigl(f - \frac{w_h - u_h^{n-1}}{\Delta t_n}\Bigr) \\
    & + ( w_h - u_h^{n-1} , w_h )  + \lambda \|w_h\|^2_{L^2(\partial\Omega)}, \nonumber
\end{align}
where we have once again added a penalization term; this time to weakly impose homogeneous Dirichlet boundary conditions on $w_h$.

Obtaining the minimizer of \eqref{eq:Phi_n} requires a maximization to evaluate \linebreak
$\phi^*_h(f - (w_h - u_h^{n-1})/{\Delta t_n})$, see \eqref{eq:phi_h}. For the remainder of this section we focus on an algorithm for solving this min-max problem. In the subsequent text it will be convenient to refer to
\begin{equation*} 
    \widetilde{\phi}_h^*(w_h; v_h) = \frac{1}{2 \kappa} \left[ \frac{ (  w_h , v_h ) }{ ( \| \nabla v_h  \|_{ L^2}^2 + \lambda \|v_h\|^2_{L^2(\partial\Omega)} )^{\frac12}} \right]^2,
\end{equation*}
which is equal to $\phi^*_h(w_h)$ if $v_h$ is a maximizer. Similarly, we write
\begin{equation} \label{eq:tilde_Phi_n}
    \widetilde{\Phi}_n(w_h; p_h) = \frac{\kappa\Delta t_n}2 \|\nabla w_h\|^2_{L^2} + \Delta t_n \, p_h + ( w_h - u_h^{n-1} , w_h )  + \lambda \|w_h\|^2_{L^2(\partial\Omega)}, 
\end{equation}
which equals $\Phi_n(w_h)$ upon choosing $p_h = \phi^*_h\bigl(f - \frac{w_h - u_h^{n-1}}{\Delta t_n}\bigr)$.

\begin{algorithm}
\caption{} \label{alg}
\begin{algorithmic}[1]
\State Compute an approximation $u^0_h \in \mathbb U_h$ to $u_0$
\For{$n=1,\ldots,N$}
\State $u^{n,0}_h \gets u^{n-1}_h$, $k \gets 0$
\While{termination\_criterion($u^{n,k}_h$, $k$) = FALSE}
\State $k \gets k+1$
\State $p_h^{n,k} = \max_{v_h \in \mathbb V_h \setminus \{0\}} \widetilde{\phi}_h^*(f - (u_h^{n,k-1} - u_h^{n-1})/\Delta t_n; v_h)$
\State $u^{n,k}_h \in \arg\min_{w_h \in \mathbb U_h} \widetilde{\Phi}_n(w_h; p_h^{n,k})$
\EndWhile
\State $u^n_h \gets u^{n, k}_h$
\EndFor
\end{algorithmic}
\end{algorithm}

We shall base the minimization of \eqref{eq:Phi_n} on Algorithm \ref{alg}. The approximation of the initial conditions in line 1 of the algorithm is a supervised learning problem. Lines 2 and 10 frame the iteration over the time steps. Lines 4 and 8 implement a loop where the optimization of $\widetilde{\phi}_h^*$ (line 6) and $\widetilde{\Phi}_n$ (line 7) are alternated.

\begin{remark}
Alternatively to the above, one could use $\phi^*(w) = \tfrac1{2\kappa} \| \nabla \Delta^{-1} w \|_{L^2}^2$ for the formulation of the method, see \cite[Example 8.104]{Roubicek13}. In this scenario we envisage $\Delta^{-1} u$ being approximated by a neural net, using an existing methodology for solving the Laplace problem.
\end{remark}

\section{BENNO: Brezis--Ekeland Neural Network Optimizer}

We make a full Python implementation \cite{github} of our deep learning approach available on Github, which is called {\em Brezis--Ekeland Neural Network Optimizer}, in short BENNO.

\subsection*{Neural Network Structure}

We describe the internal structure of the neural networks $\hat u_h$ and $\hat v_h$, which  were introduced in the previous section. Both these neural networks use five densely-connected layers with a linear activation function for the input and output layers and a (leaky) rectified linear unit activation function $\sigma(s) = \max\{0,s\} + \mu \min\{0,s\}$, $\mu \geq 0$, for the inner layers. Each layer is made of $m$ nodes, except for the output layer that presents a single node. This means that the neural network $\hat u_h$ has the following architecture:
\begin{align*} 
    S^1(x) & = W^1 x + b^1, \
    S^2(x) = \sigma( W^2 S^1(x) + b^2), \
    S^3(x) = \sigma( W^3 S^2(x) + b^3), \nonumber \\
    S^4(x) & = \sigma( W^4 S^3(x) + b^4), \
    S^5(x) = W^5 S^4(x) + b^5, \
    \text{and} \quad \hat u_h(x; \theta) = S^5(x),
\end{align*}
where the set of parameters of the neural network $\hat u_h(\cdot; \theta)$ are given by 
\[
    \theta = \{ W^1, b^1, W^2, b^2, W^3, b^3, W^4, b^4, W^5, b^5\},
\]
with $W^1 \in \mathbb{R}^{m \times d}$, $W^i \in \mathbb{R}^{m \times m}$ for $i=2,3,4 $, $W^5 \in \mathbb{R}^{1 \times m }$,  $b^j \in \mathbb{R}^m$ for $j=1,\ldots,4$ and $b^5 \in \mathbb{R}$. With a slight abuse of notation, the application of the activation function is understood elementwise: $\sigma(z)$ is the vector $(\sigma(z_1), \ldots, \sigma(z_m))$ for $z = (z_1, \ldots, z_m) \in \mathbb{R}^m$. By default, we set $\mu = 0.03$ in the definition of $\sigma$.

Also the network $\hat v_h$ has an architecture of this type; however, generally with a different parameter $m$. As indicated in the previous section, we denote the weights of $\hat v_h$ by $\eta$.

\subsection*{Adam Optimizer, Loss Functions and Algorithm}

We implemented the neural networks with the help of the Tensorflow Sequential API. Both neural networks were trained with the Adam Optimizer, a variant of the stochastic gradient descent method based on an adaptive estimation of first-order and second-order moments that improves the speed of convergence \cite{KingmaB15}. For the optimization parameters, we use the standard values $\beta_1 =0.9$, $\beta_2 = 0.999$ and $\varepsilon = 10^{-8}$ in the notation of \cite{KingmaB15}.

There are three distinct optimization scenarios with their respective loss functions:
\begin{enumerate}
    \item The approximation of $u_0$ by $u_h^0 \in \mathbb U_h$ in line 1 of Algorithm \ref{alg}: It is a supervised learning problem with the loss function $\mathcal{L} = \| u_0 -w_h  \|^2 _{L^2}$. We use the constant learning rate $\alpha = 10^{-3}$.
    \item The maximization of $\widetilde{\phi}_h^*$ in line 6 of Algorithm \ref{alg}: We use the constant learning rate $\alpha = 10^{-5}$. By default the training extends over 500 epochs.
    \item The minimization of $\widetilde{\Phi}_n$ in line 7 of Algorithm \ref{alg}: We employ the $k$-dependent decaying learning rate
    \begin{align*}
        \hspace{15mm} \alpha(k) = \, & 10^{-5} \textbf{1}_{\{k \leq 5 \}} +  10^{-6} \textbf{1}_{\{ 5 < k  \leq 50\}}+ 10^{-7} \textbf{1}_{\{ 50 < k \leq 120\}}\\
        & + 10^{-8} \textbf{1}_{\{120 < k \leq 140\}} + 10^{-9} \textbf{1}_{\{140 < k \leq 180 \}} + 10^{-10} \textbf{1}_{\{ 180 < k \}}. 
    \end{align*}
    By default the training extends over 50 epochs.
\end{enumerate}
The integrals appearing in these loss functions are evaluated with the help of a Monte-Carlo integration method, using the sampling points $\{x_i : i = 1,\ldots,N_s\} \subset \overline\Omega$. Here $N_s = N_i + N_b $, with $N_i$ points drawn from a uniform distribution in $\Omega$ and $N_b$ points drawn from a uniform distribution on $\partial\Omega$.

Finally, the default termination criterion in line 4 of Algorithm \ref{alg} is
\[
\text{termination\_criterion}(u^{n,k}_h, k) = \begin{cases}
    \text{TRUE} & : k > 200,\\ \text{FALSE} & : k \leq 200,
\end{cases}
\]
which is employed in all numerical experiments of the forthcoming section.

\section{Numerical results} \label{sec:nr}

We consider problem \eqref{eq:PDEpbl}, with $f=0$, on the domain $\Omega=(0,\pi)^d$, for $d=2,3,5, 7$. Given the initial condition $u_0(x) = \prod_{i=1}^d \sin(a_i x_i)$, for $a \in \mathbb{N}^d$ and $x \in \overline\Omega$, the exact solution to \eqref{eq:BEPheat} with $\kappa=[\sum_{i=1}^d a_i^2]^{-1}$ is $u(t,x) =  e^{-t}\prod_{i=1}^d \sin(a_i x_i)$.

We investigate the following types of approximation errors:
\begin{align*}
  MSE &= \frac{1}{N_s}\sum_{i = 1} ^{N_s} (u(t_n, x_i) - u_h^n(x_i))^2, \\
 \varepsilon_{{abs,L^ \infty}} &= \max_{i=1,\ldots,N_s}| u(t_n, x_i) - u_h^n(x_i) |, \ 
       \varepsilon_{rel, L^{2}} = \left[\frac{\sum_{i = 1} ^{N_s}( u(t_n, x_i) - u_h^n(x_i))^2} {\sum_{i = 1} ^{N_s} (u(t_n, x_i))^2}\right]^\frac12, 
\end{align*}
where MSE stands for mean square error, and the other two quantities define approximations of the $L^\infty$-norm error and of the relative $L^2$-norm error, respectively. 

In addition the Brezis--Ekeland functional itself represents a measure of the accuracy of the deep learning algorithm since we look for $u$ such that $\Phi(u) = \min \Phi = 0$. It follows that values of the loss function $\widetilde\Phi_{n}$ give us information about the quality of the training and the approximate solution $u_h^n$.

Unless otherwise stated, we use the layer width $m=m_{v_h}=30$ for the neural networks $\hat v_h$, while the the layer width $m=m_{u_h}$ will be varied for $\hat u_h$ depending on the dimension $d$. For the time discretization we use uniform time steps $\Delta t_n = \Delta t$, $n=1,\ldots, N$, where we always choose $\Delta t = 10^{-4}$. Finally, for the boundary value penalty parameter we always use $\lambda=100$.

\subsection*{Energy landscapes for a 5D problem}
Let $d=5$ and $a = (2, 2, 1, 2, 3)^{\intercal}$. We choose $N_i = 10^5$ inner and $N_b = 10^3$ boundary sampling points. The optimization to obtain the initial value approximation $u_h^0$ is done over $5\cdot10^4$ epochs. We use $m_{u_h} = 60$ for $\hat u_h$.

We are interested in the shape of the graphs of the two objective functions $\widetilde\Phi_{n}$ and $\widetilde{\phi}_h^*$ as functions of the neural network weights $\theta$ and $\eta$, respectively. This will allow us to gain insight into how challenging the training of the neural nets is. 

\begin{figure}[t]
    \begin{subfigure}[b]{0.3\textwidth}
        \centering
        \includegraphics[width = \textwidth]{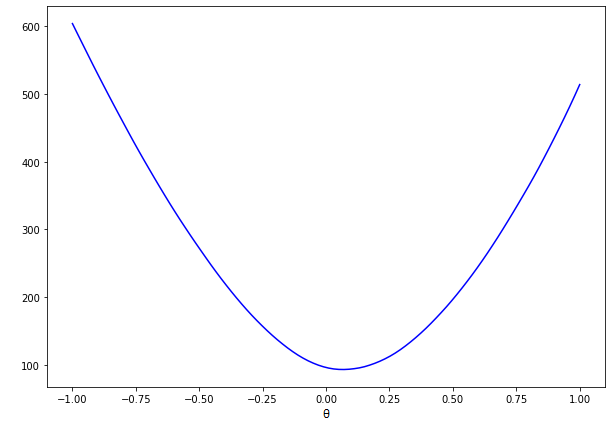}
    \end{subfigure}
    \hspace{2mm}
    \begin{subfigure}[b]{0.3\textwidth}
        \centering
        \includegraphics[width = \textwidth]{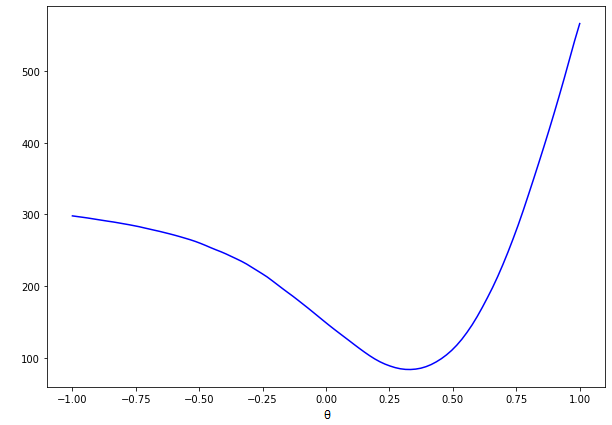}
    \end{subfigure}
    \hspace{2mm}
    \begin{subfigure}[b]{0.3\textwidth}
        \centering
        \includegraphics[width = \textwidth]{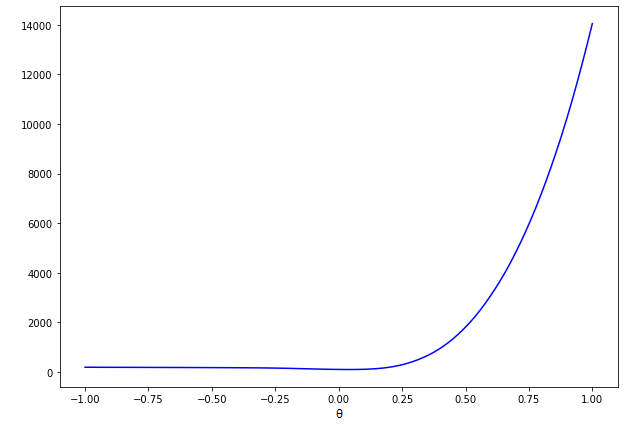}
        \end{subfigure}
        \caption{Plots of $\widetilde\Phi_{n}(\hat u_h(\theta); p_h^{1,200})$ against different components of $\theta$: $b^1_{60}$ (left), $W^3_{46,60}$ (middle) and $W^3_{7, 45}$ (right).}
        \label{BE_landscape}
\end{figure}

Having computed $u_h^1 = u_h^{1,K}=\hat u_h(\hat\theta)$, in Figure~\ref{BE_landscape} we plot the loss function $\widetilde\Phi_{n}(\hat u_h(\theta); p_h^{n,K})$, for $n=1$ and $K=200$, against selected entries of $\theta$. In particular, for each plot we keep all the weights in $\theta=\hat\theta$ fixed, apart from a single entry of $\theta$, that we continuously vary from $-1$ to $1$. In this way it is possible to visualize how the Brezis--Ekeland functional varies depending on certain parameters of the neural network $\hat u_h$. While generally smooth, we note that the right plot in Figure~\ref{BE_landscape} shows that $\widetilde\Phi_n$ has a nearly vanishing gradient when the parameter $W^3_{7,45}$ varies in $[-1,0]$, which may require attention during the optimization process.

\begin{figure}[t]
    \begin{subfigure}[b]{0.3\textwidth}
        \centering
        \includegraphics[width = \textwidth]{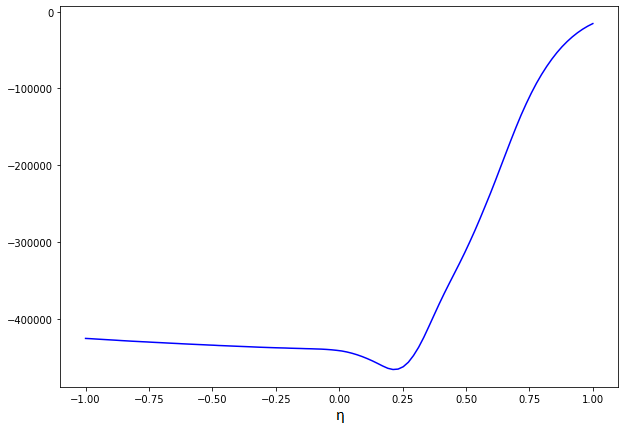}
    \end{subfigure}
    \hspace{2mm}
    \begin{subfigure}[b]{0.3\textwidth}
        \centering
        \includegraphics[width = \textwidth]{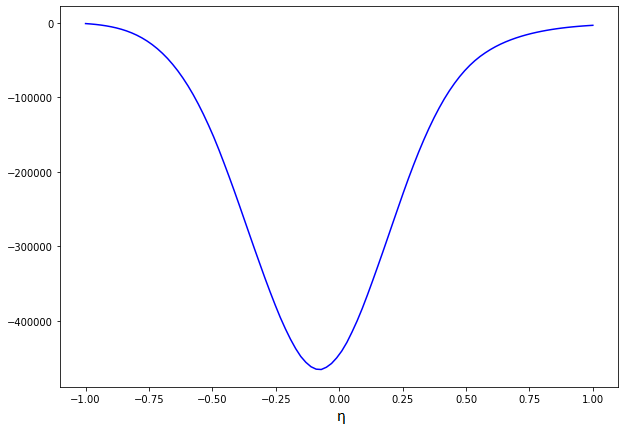}
    \end{subfigure}
    \hspace{2mm}
    \begin{subfigure}[b]{0.3\textwidth}
        \centering
        \includegraphics[width = \textwidth]{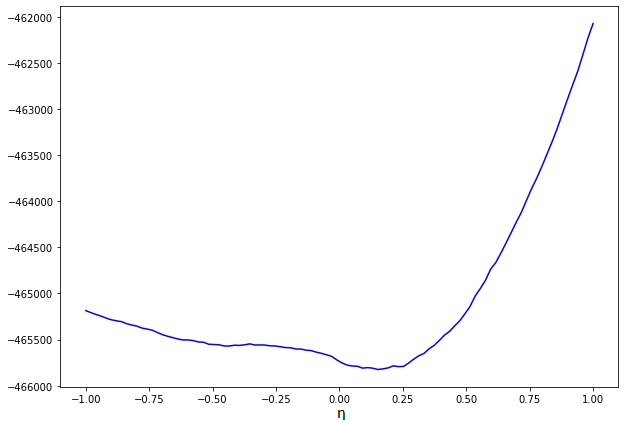}
        \end{subfigure}
        \caption{Plots of $  - \widetilde{\phi}_h^*( (u_h^{1,199}-u_h^{0})/\Delta t; \hat v_h(\eta)) $ against different components of $\eta$: $b^2_{30}$ (left),  $W^3_{2,1}$ (middle) and $W^4_{21,3}$ (right).}
        \label{negNorm landscape}
\end{figure}

Similarly, in Figure~\ref{negNorm landscape} we show the loss function $-\widetilde{\phi}_h^*( (u_h^{n,k-1}-u_h^{n-1})/\Delta t; \hat v_h(\eta)) $ plotted against selected entries of the neural network weights $\eta$. Once again we observe nearly flat parts in the graph, but now in addition we see also some non-convex and non-smooth regions, which may pose challenges during the optimization.

\subsection*{Error quantities for a 5D problem}
We use the previous example to compute error quantities for the trained neural networks. Table~\ref{Test5.1} shows values of the Brezis--Ekeland loss function ($\widetilde\Phi_{n}$), the mean square error (MSE), the absolute error ($\varepsilon_{{abs,L^ \infty}}$) and the relative error ($\varepsilon_{rel, L^{2}}$) for the neural networks $u_h^n$ for every $t_n$, $n = 0,\ldots, N$. The reported values of $\widetilde\Phi_{n}$ are large, possibly caused by the fact that the measure of $\Omega$ is considerable with $|\Omega|=\pi^5 \approx 306$, as is the scaling of the boundary term $\lambda \| \cdot \|_{L^2(\partial \Omega)}^2$ with $\lambda = 100$ and $|\partial \Omega|= 2 \cdot 5 \cdot \pi^4 \approx 974$. We observe that the MSE and $\varepsilon_{rel, L^{2}}$ increase with time, which is typical for approximations of parabolic PDEs. The reduction of $\varepsilon_{{abs,L^\infty}}$ in time can be explained by the $e^{-t}$ scaling of the true solution.
\begin{table}[t]
    \small \centering
    \begin{tabular}{||c||cccc||}
    \hline
        & $\widetilde\Phi_{n}$  & MSE & $\varepsilon_{{abs,L^ \infty}}$ & $\varepsilon_{rel, L^2}$\\ [0.5ex] \hline
        $t_0$  & ---    & 5.154e-04 & 0.128 & 0.129\\
        $t_1$  & 99.450 & 5.460e-04 & 0.125 & 0.133\\
        $t_2$  & 62.899 & 5.779e-04 & 0.128 & 0.136\\
        $t_3$  & 62.221 & 6.235e-04 & 0.128 & 0.140\\
        $t_4$  & 53.464 & 6.608e-04 & 0.149 & 0.146\\
        $t_5$  & 62.384 & 7.032e-04 & 0.152 & 0.150\\
        $t_6$  & 49.826 & 7.418e-04 & 0.162 & 0.155\\
        $t_7$  & 47.909 & 7.742e-04 & 0.164 & 0.158\\
        $t_8$  & 44.545 & 8.020e-04 & 0.172 & 0.161\\
        $t_9$  & 45.276 & 8.351e-04 & 0.174 & 0.164\\
        $t_{10}$ & 43.170 & 8.468e-04 & 0.177 & 0.165\\
        \hline
    \end{tabular}
    \caption{Error quantities at times $t_n$, $n = 0,\ldots, N$, for the 5D test problem.}
    \label{Test5.1}
\end{table}

\begin{figure}[t]
    \begin{subfigure}[b]{0.48\textwidth}
        \centering
        \includegraphics[width = \textwidth]{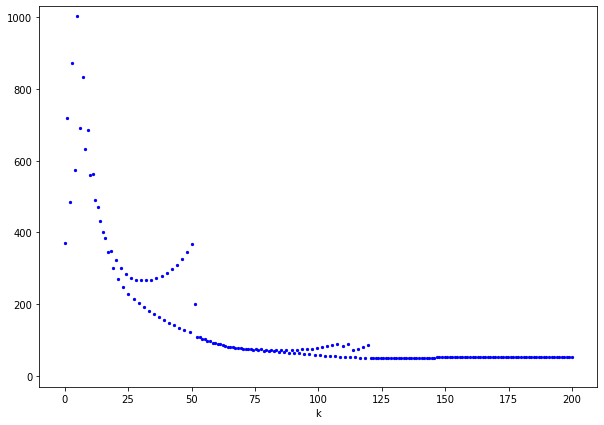}
    \end{subfigure}
    \hspace{1mm}
    \begin{subfigure}[b]{0.48\textwidth}
        \centering
        \includegraphics[width = \textwidth]{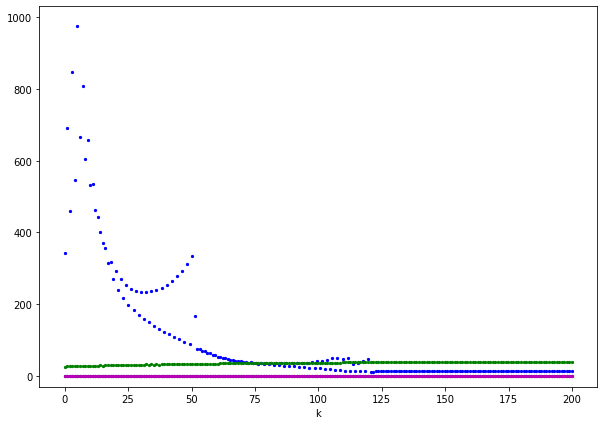}
    \end{subfigure}
    \caption{Plots of $\widetilde\Phi_{n}(u_h^{n,k}; p_h^{n,k})$ (left) and of the 4 terms contributing to it (right) against $k$, for $n=4$, for the 5D test problem. For the right plot the chosen colours are red for $ \frac{\kappa\Delta t}2 \|\nabla u_h^{n,k}\|^2_{L^2}$, blue for $ \Delta t\, p_h^{n,k}$, magenta for $( u_h^{n,k} - u_h^{n-1} , u_h^{n,k} )$ and green for $\lambda \|u_h^{n,k}\|^2_{L^2(\partial\Omega)}$.}
    \label{Be and NegNorm for 5th time iter}
\end{figure}

Apart from the global error properties, we are also interested in how these quantities change during the training process. In Figure~\ref{Be and NegNorm for 5th time iter} we plot the loss function $\widetilde\Phi_{n}(u_h^{n,k}; p_h^{n,k})$ and the four contributions to it against $k$, for $k=1, \ldots, K = 200$, during the training for the time $t_4 = 4\cdot10^{-4}$, i.e.\ $n=4$. We observe a significant  decrease of the Brezis--Ekeland functional $\widetilde\Phi_{n}$ during the training, from about 400 for $k=1$ to about $50$  for $k=200$, when the weights seem to have converged. Observe also that the decrease is non-monotone, with a global maximum of about $1000$, and that the graph is rather oscillatory. In addition, we note that after an initial increase, the functional decays rapidly at first and then slower as the iteration proceeds. The plot of the four contributions reveals that the term $\lambda \|u_h^{n,k}\|^2_{L^2(\partial\Omega)}$ is the dominant contribution in the Brezis--Ekeland functional \eqref{eq:tilde_Phi_n} once the iterative scheme settles down.

\begin{figure}[t]
    \begin{subfigure}[b]{0.48\textwidth}
        \centering
        \includegraphics[width = \textwidth]{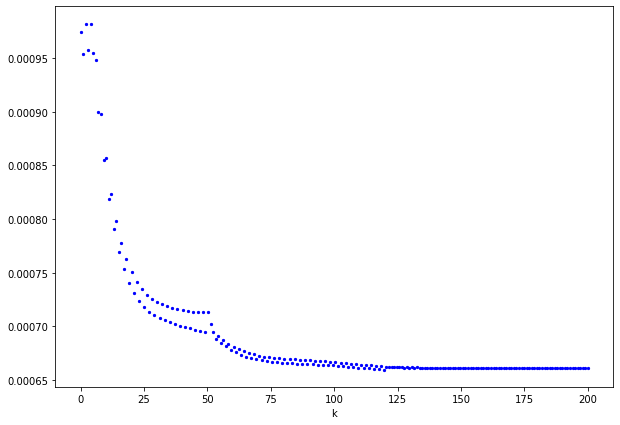}
    \end{subfigure}
    \hspace{1mm}
    \begin{subfigure}[b]{0.48\textwidth}
        \centering
        \includegraphics[width = \textwidth]{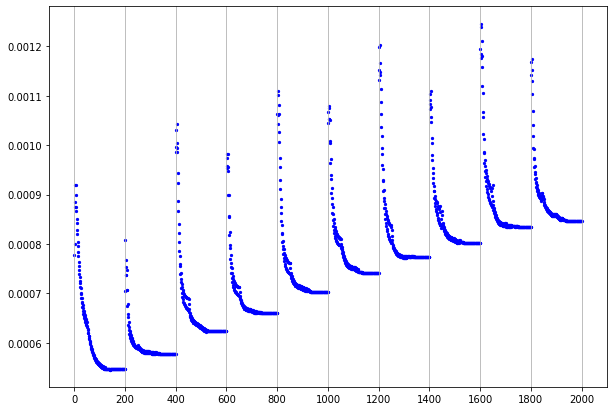}
    \end{subfigure}
    \caption{Plots of the MSE against $k$ at time $t_4 = 4\cdot10^{-4}$ (left), and of the MSE against $(n-1)K+k$, for 
    $n=1,\ldots,N$. (right)}
    \label{MSE for 5th and total time iter}
\end{figure}

The analogous plot for the mean square error (MSE) is shown on the left of Figure~\ref{MSE for 5th and total time iter}, where we again notice an oscillatory decrease until convergence is reached. In addition, on the right of Figure~\ref{MSE for 5th and total time iter} we show the concatenated plots of the MSE against $k$, for every time $t_n$, $n=1, \ldots, N$. To help differentiate the different time steps, we indicate the start of the training for a new time step with vertical lines. The figure demonstrates that overall the MSE increases in time, but that each training procedure decreases the MSE until convergence can be observed.

\subsection*{Dependence on the dimension $d$}
Here we let $d=2$, $3$ or $5$, and set $a= (2, 2)^\intercal$, $a = (2, 2, 3)^\intercal$ and $a = (2, 2, 1, 2, 3)^\intercal$, respectively. The number of sampling points in $\Omega$ is $N_i = 10^4$ for $d=2$ and $N_i = 10^5$ for $d=3,5$, while $N_b = 400, 600, 1000$ for $d=2,3,5$, respectively. Moreover, for the training of the initial conditions $u_h^0$ we use $5\cdot10^3$ epochs in the cases $d=2,3$, and $5\cdot10^4$ for $d=5$. We set $m_{u_h}= 60$ throughout.

Table~\ref{Table: Comparison between 2D, 3D and 5D tests.}  shows the values of the Brezis--Ekeland loss function $\widetilde\Phi_{n}$, the MSE, the absolute error $\varepsilon_{{abs,L^ \infty}}$, the relative error $\varepsilon_{rel, L^{2}}$ and the GPU time used for the training of the neural networks $\hat u_h$ and $\hat v_h$, for the tenth time step, $t_{10} = 0.001$, for the three problems $d=2,3,5$. It can be observed that an increase in the dimension leads to a decrease in the accuracy of the algorithm. Indeed, reading the table from left to right, we notice that all the measures increase monotonically with the dimension, with the only exception being the absolute error $\varepsilon_{{abs,L^ \infty}}$, which actually reduces slightly going from the 3D to the 5D problem. While a growth of the error is to be expected with an increase of the dimension, the figures here indicate that further improvements of the neural net architecture and the training methodology should be investigated.

\begin{table}[t]
    \small \centering
    \begin{tabular}{||c||c|c|c||}
    \hline
        & 2D & 3D & 5D \\ \hline
        $\widetilde\Phi_{n}$ & 4.54e-02 & 3.97 & 39.93 \\
        MSE & 1.74e-04 & 1.00e-03 & 7.67e-03 \\ 
        $\varepsilon_{{abs,L^ \infty}}$  & 6.25e-03 & 0.14 & 0.13 \\ 
        $\varepsilon_{rel, L^{2}}$  & 2.71e-02 & 8.98e-02 & 0.15\\
         GPU time [s] & 58195 & 73595 & 100641\\
        \hline
    \end{tabular}
    \caption{Comparison of error quantities and the GPU time for the three test problems with $d=2,3,5$.}
    \label{Table: Comparison between 2D, 3D and 5D tests.}
\end{table}

\subsection*{Effect of the number of nodes $m$ for a 7D problem}
Let $d=7$ and $a = (2, 2, 1, 3, 2, 2, 3)^\intercal$. We choose $N_i = 10^5$ and $N_b = 1400$ sample points. The initial conditions are training over $5\cdot10^4$ epochs. The layer width $m_{u_h}$ of the networks $\hat u_h$ varies between $60$ and $100$.

From the error quantities reported in Table~\ref{Table: Comparison between 7D with 60 and 70 neurons.} we can immediately see that using more nodes in the network architecture is beneficial in higher dimensions. In fact, all the reported quantities are lower when considering the case with $m_{u_h} = 100$. 

\begin{table}[t]
    \centering
    \begin{tabular}{||c||c|c||}
    \hline
        $m_{u_h}$ & 60 & 100 \\
        \hline  
        $\widetilde\Phi_{n}$ & 469.11 & 60.44 \\
        MSE & 3.80e-03 & 1.21e-03 \\ 
        $\varepsilon_{{abs,L^ \infty}}$  & 0.59 & 0.31 \\ 
        $\varepsilon_{rel, L^{2}}$  & 0.70 & 0.38 \\
        GPU time [s] & 119972 & 117657 \\
        \hline
    \end{tabular}
    \caption{Comparison of error quantities and the GPU time for the 7D test problem with either 60 or 100 nodes per layer.}
    \label{Table: Comparison between 7D with 60 and 70 neurons.}
\end{table}

\begin{figure}[!ht]
    \begin{subfigure}[b]{0.48\textwidth}
        \centering
        \includegraphics[width = \textwidth]{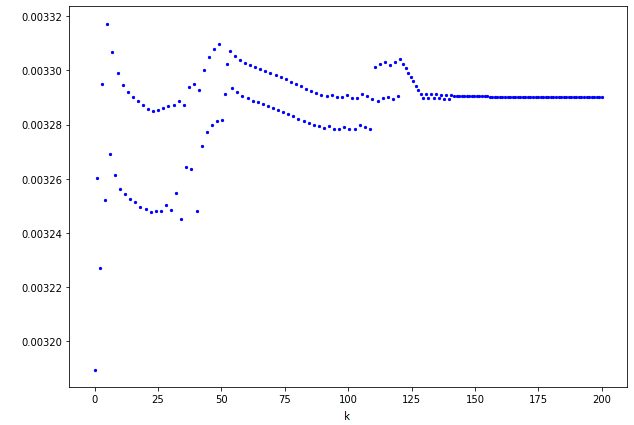}
    \end{subfigure}
    \hspace{1mm}
    \begin{subfigure}[b]{0.48\textwidth}
        \centering
        \includegraphics[width = \textwidth]{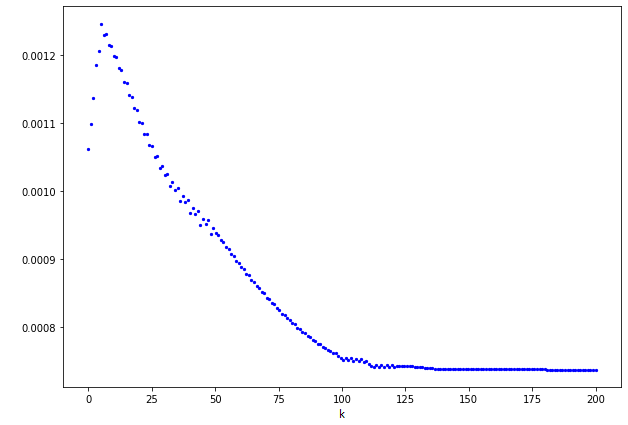}
    \end{subfigure}
    \caption{Plots of the MSE against $k$ at time $t_3=3\cdot10^{-4}$ for the 7D problem with neural netwoks using 60 (left) and 100 (right) nodes per layer.}
    \label{Fig:Comparison between 7D tests}
\end{figure}

Finally, in Figure~\ref{Fig:Comparison between 7D tests} we present plots of the MSE against $k=1, \ldots, K = 200$ for the training of the third time step, $n=3$, comparing the performance of the two different network structures. We note that the MSE is on average increasing for the network with 60 nodes per layer, possibly indicating that the neural network does not have enough expressivity to make the learning effective. In contrast, when training with wider layers one is able to decrease the MSE after the typical initial increase. These results strongly indicate that for higher dimensional problems, more elaborate networks perform better in practice. Unfortunately, memory and GPU time limitations mean that at present we are not able to investigate this trend for even higher dimensional problems.

\section{Conclusions}
We introduced a novel deep learning approach for the
numerical solution of PDEs using the Brezis--Ekeland principle. As a proof of concept we implemented a practical algorithm for the heat equation and presented results for experiments up to dimension 7. Higher dimensional problems are particularly computationally challenging, and more research into the optimal design for the employed neural networks is needed. In addition, an extension of the implemented method to nonlinear problems is part of future research.

\bibliographystyle{equadiff}
\bibliography{refs.bib}
\end{document}